\documentclass[a4paper]{amsart}
\usepackage{a4wide}
\usepackage{natbib}
\newtheoremstyle{mathfin}{3pt}{3pt}{\itshape}{}{\sc}{.}{.5em}{}
\theoremstyle{mathfin}
\newtheorem{lemma}{Lemma}[section]
\newtheorem{theorem}{Theorem}[section]
\newtheorem{example}{Example}[section]
\begin{document}
\title[A general proof of the Dybvig-Ingersoll-Ross-Theorem]{A general proof
of the Dybvig-Ingersoll-Ross-Theorem: Long forward rates can never
fall}
\author{Friedrich Hubalek}
\author{Irene Klein}
\author{Josef Teichmann}

\thanks{We thank Walter Schachermayer for drawing our attention
to the topic of this note and Freddy Delbaen for
his comments and improvements in general, and especially
regarding the statement
and proof of Lemma~\ref{Freddy}.}

\address{Friedrich Hubalek,
         Department of Financial and Actuarial Mathematics,
         Vienna University of Technology,
         Wiedner Hauptstrasse 8--10, 1040 Vienna, Austria.}
\email{fhubalek@fam.tuwien.ac.at}

\address{Irene Klein,
´        Department of Statistics and Decision Support Systems,
         University of Vienna,
         Br\"unnerstrasse 72, 1210 Vienna, Austria.}
\email{irene.klein@univie.ac.at}

\address{Josef Teichmann,
         Department of Financial and Actuarial Mathematics,
         Vienna University of Technology,
         Wiedner Hauptstrasse 8--10, 1040 Vienna, Austria.}
\email{josef.teichmann@fam.tuwien.ac.at}



\begin{abstract}
A general proof of the Dybvig-Ingersoll-Ross Theorem on the
monotonicity of long forward rates is presented. Some inconsistencies
in the original proof of this theorem are discussed.
\end{abstract}

\maketitle

\section{Introduction}

It is an interesting question to analyse the stochastic nature of long
term rates in interest rate markets. In \citet*{DIR:96} the authors
show that long forward and zero coupon rates can never fall. In their
proof they implicitly use an ``ergodicity'' assumption, which is
economically reasonable, but does not hold in any arbitrage-free
interest rate model (see Example~4.1). Furthermore there are some
difficulties with a limiting procedure that are addressed in
\citet*{McC:00}, but even under that ``ergodicity'' assumption the
strategy in the proof of \citet*{McC:00} is anticipative, so not
admissible for a no-arbitrage argument. In this note we prove  without
any additional assumption that long forward rates  can
never fall, if they exist.

\section{Interest Rate models}
Suppose we are given a probability
space $(\Omega,\mathcal{F},P)$ with
filtration $(\mathcal{F}_{t})_{t\geq0}$ where the time parameter is either
discrete ($t\in\mathbb{N}$) or continuous ($t\in\mathbb{R}_{\geq0}$).
Prices of default-free
zero coupon bonds $P(t,T)$ are modelled as semimartingales for $0\leq t\leq T$
with respect to $(\mathcal{F}_{t})_{0\leq t\leq T}$. The process
$\{P(t,T)\}_{0\leq t\leq T}$ is strictly positive, furthermore we assume the
normalization $P(T,T)=1$. No arbitrage in this setting is usually given by the
following requirement, which we shall assume throughout: there exists a
probability measure $Q$ and the $(\mathcal{F}_{t})_{t\geq0}$-adapted interest
rate process $(R_{t})_{t\geq0}$ (the rates can be negative, too) such that the
following conditions hold:
\begin{enumerate}
\item $B_{t}$ is a well-defined predictable, strictly positive process.
In the discrete case it is defined by
\[
B_{t}:=\prod_{i=0}^{t-1}(1+R_{i}),
\]
in the continuous case by
\[
B_{t}:=\exp\left(\int_{0}^{t}R_{s}ds\right).
\]
Furthermore we assume that $B_t/B_{t+h}$ are integrable with respect
to $Q$ for $t\geq0$ and $h\geq0$.

\item  The measure $Q$ is locally equivalent to $P$, i.e. $Q|\mathcal{F}%
_{t}\sim P|\mathcal{F}_{t}$ for $t\geq0$. Remark that we did not assume
usual conditions for the filtration, see \citet*{Del:93} for a
related discussion.

\item  The discounted processes $B_{t}^{-1}P(t,T)$ are $Q$-martingales for
$0\leq t\leq T$.
\end{enumerate}

The no arbitrage condition yields therefore%
\[
B_{s}^{-1}P(s,T)=E_Q(B_{t}^{-1}P(t,T)|\mathcal{F}_{s})
\]
for $t\geq s$. From the given normalization and adaptedness
of $B_{t}$ the representation%
\[
P(t,T)=E_Q\left(\frac{B_{t}}{B_{T}}\bigg|\mathcal{F}_{t}\right)
\]
of the price processes follows.

\subsection{Discrete Case}

The forward rate process $f(t,T)$ is well-defined by the following formula%
\[
P(t,T)=(1+f(t,T))P(t,T+1)
\]
for $0\leq t\leq T$, the zero coupon rate $z(t,T)$ is given by the formula%
\[
P(t,T)=\frac{1}{(1+z(t,T))^{T-t}}%
\]
Both processes are $(\mathcal{F}_{t})_{t\geq0}$-adapted, however, their
integrability properties depend on the price processes. We obtain furthermore
the identification $R_{t}=f(t,t)$.

\begin{lemma}
Assume that the long forward rate exists as almost sure limit, i.e.%
\[
\lim_{T\rightarrow\infty}f(t,T)=f_{L}(t)
\]
then the long zero coupon rate%
\[
\lim_{T\rightarrow\infty}z(t,T)=z_{L}(t)
\]
exists as almost sure limit and $z_{L}(t)=f_{L}(t)$.
\end{lemma}

\begin{proof}
The proof can be
found in \citet*{DIR:96}, too. One applies the formula%
\[
z(t,T)=\frac{1}{P(t,T)^{\frac{1}{T-t}}}-1
\]
where we insert%
\[
P(t,T)=\prod_{I=t}^{T-1}(1+f(t,I))^{-1}%
\]
which yields the result.
\end{proof}

\subsection{Continuous Case}

The forward rate is $f(t,T)$ defined via the following formula%
\[
P(t,T)=\exp\left(-\int_{t}^{T}f(t,s)ds\right)
\]
We assume that the forward rate exists as an adapted process. We obtain
under regularity assumptions $f(t,t)=R_{t}$ for $t\geq0$. The analogue
to the zero
coupon rate is given by the yield process%
\[
z(t,T):=\frac{1}{T-t}\int_{t}^{T}f(t,s)ds.
\]

\section{Long Forward Rates never fall}

Considering a unifying approach to discrete and continuous time interest rate
models we can write in the above notions%
\[
P(t,T))=E_Q(B_{T}^{-1}B_{t}|\mathcal{F}_{t})
\]
for $t\leq T$. We assume that $ z_{L}(t) $ is an almost surely finite
random variable: the process $\{z_{L}(t)\}_{t\geq0}$ is increasing (in
the sense that $ z_{L}(t) \geq z_{L}(s) $ a.s.\ for $ t \geq s $)
if and only if%
\[
x_{L}(t)=\lim_{T\rightarrow\infty}P(t,T)^{\frac{1}{T}}%
\]
is decreasing, since in the discrete case
$z_{L}(t)=\frac{1}{x_{L}(t)}-1$ and in the continuous case
$z_{L}(t)=-\ln x_{L}(t)$. We denote by $ x(t,T) $ the random variable $
P(t,T)^{\frac{1}{T}} $. The almost sure existence of $z_{L}(t)$ is
equivalent to the existence of $x_{L}(t)$ for all $t\geq0$ by
definition.

\begin{theorem}\label{dir}

If $x_L(t)$ and $x_L(s)$ exist almost surely for $t \geq s \geq0$ then
$$ x_L(s)\geq x_L(t)~\mbox{a.s.} $$
\end{theorem}

For the proof of this theorem we apply the following technical lemma,
which generalizes the well-known
fact $\lim_{p\to\infty}
\left\|X\right\|_p=\left\|X\right\|_\infty$ for $X\in L^\infty$.
\begin{lemma}\label{Freddy}
Let $\{X_n\}_{n\geq0}$ be a sequence of non-negative random variables
on a probability space $(\Omega,\mathcal F,P)$ and let $\mathcal G$ be
a sub-$\sigma$-algebra of $\mathcal F$. Suppose $X_n$ converges to the
random variable $X$ a.s.\ and
$
\liminf_{n\to\infty}E[X_n^n|\mathcal G]^{\frac1n}=C<\infty~\mbox{a.s.}
$
Then $X\leq C$ a.s.
\end{lemma}
\begin{proof}
Replacing $X$ by $X1_A$, where $A=\{C\leq k\}\in\mathcal G$, letting
$k\to\infty$ allows us to replace~$C$ by a bounded $\mathcal
G$-measurable random variable. Take $f\geq0$, bounded, $E[f]=1$. Using
the conditional Fatou Lemma for a.s.\ convergence and the conditional
H\"older Inequality we obtain
\begin{eqnarray*}
E[Xf]&=&E[\lim X_nf]=E[E[\lim X_nf|\mathcal G]] \leq E[\liminf
E[X_nf|\mathcal G]]\\ &\leq&E\left[\liminf E[X_n^n|\mathcal
G]^{\frac1n} E\left[f^{\frac{n}{n-1}}|\mathcal
G\right]^{\frac{n-1}{n}}\right] \leq E[CE[f|\mathcal G]],
\end{eqnarray*}
since $$\lim_{n\to\infty} E\left[f^{\frac{n}{n-1}}|\mathcal
G\right]^{\frac{n-1}{n}} =E[f|\mathcal G] $$ by the conditional
Lebesgue Dominated Convergence Theorem. Finally $E[CE[f|\mathcal G] = E[Cf]$,
since $C$ is $\mathcal G$-measurable, hence $X\leq C$, since $f$ was
arbitrary $\mathcal F$-measurable.
\end{proof}

\begin{proof}
Now we can prove Theorem \ref{dir}. Therefore we fix $t\geq s$, we have
to prove
\[
\lim_{T\rightarrow\infty}P(s,T)^{\frac{1}{T}}\geq\lim_{T\rightarrow\infty
}P(t,T)^{\frac{1}{T}}.
\]
By conditioning we have
\[
P(s,T)=E\left(  \left.  \frac{B_{s}}{B_{t}}E\left(  \left.  \frac{B_{t}}%
{B_{T}}\right|  \mathcal{F}_{t}\right)  \right|  \mathcal{F}_{s}\right)
=E\left(  \left.  \frac{B_{s}}{B_{t}}P(t,T)\right|
\mathcal{F}_{s}\right).
\]
We define $\widetilde{Q}$ by
\[
\frac{d\widetilde{Q}}{dQ}=\frac{1}{P(s,t)}\frac{B_{s}}{B_{t}}.%
\]
This measure is the forward (time~$s$) neutral measure for
maturity~$t$. We can write
\[
\frac{P(s,T)}{P(s,t)}=\tilde{E}\left(  \left.  P(t,T)\right|  \mathcal{F}%
_{s}\right),
\]
with $\tilde E$ denoting expectation with respect to~$\tilde Q$.
We have
\[
x_{L}(s) = \lim_{T\rightarrow\infty}P(s,T)^{\frac{1}{T}}=\lim_{T\rightarrow\infty}%
\tilde{E}\left(  \left.  x(t,T)^{T}\right|  \mathcal{F}_{s}\right)
^{\frac
{1}{T}}%
\]
and the question reduces to
\[
\lim_{T\rightarrow\infty}\tilde{E}\left(  \left.  x(t,T)^{T}\right|
\mathcal{F}_{s}\right)
^{\frac{1}{T}}\geq\lim_{T\rightarrow\infty}x(t,T).
\]
which is a consequence of Lemma \ref{Freddy}.
\end{proof}

\section{Comments}
In \citet*{DIR:96} the statement, that the long zero coupon rate in a
discrete interest rate model can never fall, is proved by a
no-arbitrage argument. The constructed strategy is non-anticipative,
only if the zero coupon rates satisfy an additional assumption:  The
authors assume implicitly that $ z_L (t)$ is somehow $ \mathcal{F}_s
$-measurable for $t>s$. This ``ergodicity'' assumption is
economically reasonable, since the long rates should not depend on the
time point $t$ where we observe them. Nevertheless this does not hold
in all interest rate models.

We provide the following well-known
example (see \citet*{IngSkeWei:78})
to show that there exist ``non-ergodic'' interest rate models.
\begin{example}
We take
$r_{t}=r_{0}+\delta N_{t}$ where $N_{t}$ is a Poisson process with
intensity $\lambda$, jump size
$1$ and $\delta>0$ in its natural filtration. In this case%
\[
z(t,T)=r_{t}+\lambda-\frac{\lambda}{\delta(T-t)}(1-e^{-\delta(T-t)})
\]
which yields%
\[
z_{L}(t)=r_{t}+\lambda.
\]
This process is increasing, but the model is not ``ergodic'', since it
generates the filtration.
\end{example}

\end{document}